\documentclass[a4paper,12pt]{amsart}
\usepackage{amsfonts} 
\usepackage{amssymb} 
\usepackage{amsmath,amsthm} 
\usepackage{mathptm}
\usepackage[mathscr]{eucal} 
\usepackage{mathrsfs} 

\newtheorem*{burns}{Theorem}  
\newtheorem{theoreme}{Theorem}[section] 

\newtheorem{lemma}[theoreme]{Lemma} 
\newtheorem{prop}[theoreme]{Proposition} 
\newtheorem{cor}[theoreme]{Corollary} 

\theoremstyle{definition}
\newtheorem{definition}[theoreme]{Definition}

\theoremstyle{remark}
\newtheorem{remark}[theoreme]{Remark}
\newtheorem{consequence}[theoreme]{Consequence}

\def \Rk {\ {\bf Remark.} } 
 
\def \sm {\smallsetminus } 
\newcommand{\be}{\begin{enumerate}}  \newcommand{\ee}{\end{enumerate}} 
\newcommand{\bi}{\begin{itemize}}  \newcommand{\ei}{\end{itemize}} 
\newcommand{\bd}{\begin{description}}  \newcommand{\ed}{\end{description}} 
 
\newcommand{\lan}{\langle}   \newcommand{\ran}{\rangle} 
\newcommand{\norm}[1]{\lVert #1\rVert}
 
\newcommand{\comment}[1]{}

\def \R {\mathbb{R}} \def \C {\mathbb{C}}

\def \harm {\mathscr{H}} 
 
\DeclareMathOperator{\ke}{Ker}  
\DeclareMathOperator{\im}{Im} 

\numberwithin{equation}{section}  
\renewcommand{\phi}{\varphi} 
\renewcommand{\epsilon}{\varepsilon} 
 
\title{Embeddability of some strongly pseudoconvex CR manifolds} 
\author{G. Marinescu} 
\address{Humboldt-Universit{\"a}t zu Berlin, Institut f{\"u}r Mathematik, 
Unter den Linden 6,10099 Berlin, Germany} 
\email{george@mathematik.hu-berlin.de} 
\author{N. Yeganefar}
\address{Universit\'e de Nantes, D\'epartement de Math\'ematiques, 2 rue de la Houssini\`ere, 
BP 92208, 44322 Nantes cedex 03, France}
\email{nader.yeganefar@math.univ-nantes.fr}
\thanks{The second-named author has been (partially) supported by the European Commission 
through the Research Training Network  HPRN-CT-1999-00118  "Geometric Analysis".} 
\date{\today} 

\begin{document} 

\begin{abstract}
We obtain an embedding theorem for compact strongly pseudoconvex CR manifolds which are bounadries of 
some complete Hermitian manifolds. We use this to compactify some negatively curved K\"ahler manifolds 
with compact strongly pseudoconvex boundary. An embedding theorem for Sasakian manifolds is also derived.  
\end{abstract} 

\maketitle 

\section{Introduction} 

A natural question in the study of CR manifolds is the embeddability question: given a strongly pseudoconvex 
CR manifold $X$, do there exist an integer $N$ and a smooth embedding $M \hookrightarrow \C ^N$ such that the 
CR structure on $X$ is the one 
induced from the complex structure on $\C ^N$?  
It is known that the embedability of $X$ is equivalent to the fillability, 
that is the existence of a compact complex manifold $T$ with strongly  
pseudoconvex boundary, such that $\partial T=X$. Indeed, Harvey-Lawson's  
theorem \cite{HaLa:75} shows that $X$ is fillable if it is embeddable. 
Conversely, if $X$ is fillable and $\partial T=X$, a theorem of Grauert  
\cite{Gra:62} 
shows that $T$ can be taken to be a Stein space with isolated singularities, 
and one can apply the embedding theorem of Remmert-Bishop-Narasimhan  
\cite{Na:60}. 
 
By a theorem due to Boutet de Monvel  
\cite{Bo}, any compact $(2n-1)$-dimensional strongly pseudoconvex CR  
manifold is embeddable in some $\C ^N$ if its dimension is greater than or equal  
to five ($n\geqslant 3$).  
On the other hand, there is a classical example of H. Rossi \cite{Ro} 
(see also \cite{Bu}, \cite{Fa:92}) 
which shows that an arbitrarily small, real-analytic perturbation of the standard structure on  
the three-sphere may fail to be embeddable.  
There exists an extensive literature on the embeddability 
of deformations of CR structures, 
see e.g. J. Bland and C. L. Epstein \cite{Bl-E},  
D. M. Burns and C. L. Epstein \cite{Bu-E}, C. L. Epstein and G. M. Henkin \cite{E-H}, 
L. Lempert 
\cite{Le2}. 
 
Our goal is to get embedding theorems for compact strongly pseudoconvex 
CR manifolds, which are already boundaries of non-compact 
strongly pseudoconvex complex manifolds (see Definition \ref{bdy}). 
The complex manifolds we consider are complete  
Hermitian manifolds with an additional assumption on the spectral  
properties of the operator $\bar\partial {+\overline\partial }^*$ acting on  
square integrable forms.  
Our main technical result is the following theorem: 
\begin{theoreme}\label{embeddingintro} 
Let $M$ be a complex manifold with compact strongly pseudoconvex boundary 
and of complex dimension $n\geqslant 2$.  
Assume that $\operatorname{Int}M$ has a complete Hermitian metric such that  
$(\overline{\partial}+\overline{\partial}^*)^2$ is Fredholm 
in bidegree $(0,1)$. Then there exists a CR 
embedding of the boundary $\partial M$ in $\C ^N$ for some  
integer $N$. 
\end{theoreme} 
 
We then give two applications of this result. The first one concerns 
manifolds with pinched negative curvature and compact strongly pseudoconvex  
boundary. 
\begin{theoreme}\label{pinched} 
Let $M$ be a connected complex manifold with compact strongly pseudoconvex  
boundary and of complex dimension $n\geqslant 2$. Assume that $\operatorname{Int}M$ is endowed 
with a complete K\"ahler  
metric with pinched negative curvature, such that away from a neighborhood of 
$\partial M$, the volume of $M$ is finite. Then  
\begin{enumerate} 
\item $\partial M$ is embeddable in some $\C^N$. 
\item $M$ can be compactified to a strongly pseudoconvex domain in a projective variety 
by adding an exceptional analytic set, that is,  
there exists a compact strongly pseudoconvex domain $D$ in a smooth projective variety and
an embedding $h:M\to D$ which is a biholomorphism between $\operatorname{Int}M$
and $h(\operatorname{Int}M)$, $h(\partial M)=\partial D$,
and $D\smallsetminus h(M)$ is an exceptional analytic set which can be blown down to a finite set.  
\end{enumerate} 
\end{theoreme}  
This is a generalization of \cite[Main Theorem]{SY:82} in the case when a strongly
pseudoconvex end is allowed. The method follows \cite{Na-R} and consists in filling 
$\partial M$ with a Stein space $S$
(this is granted by Theorem \ref{embeddingintro}) and using a theorem of Lempert 
\cite[Theorem 1.2]{Le2} to embed $S$ in a projective variety $X$. The complement
$Y:=X\smallsetminus S$ can be glued to the original manifold $M$ to construct a 
projective pseudoconcave manifold $\widehat{M}$. By using the results of \cite{An:63} 
we can compactify this manifold and by \cite{SY:82} or \cite{M-D} the complement of 
$\widehat{M}$ is an exceptional analytic set. 

As a consequence we study some  
quotients of the unit complex ball $B$ in $\C ^n$ which where considered  
by Burns and Napier-Ramachandran \cite[Theorem 4.2]{Na-R}. 
\begin{cor}\label{Burnsintro} 
Let $\Gamma$ be a torsion-free discrete group of automorphisms of the unit  
ball $B$ in $\C ^2$ and let $M=B/\Gamma$. Assume that the limit set  
$\Lambda$ is a proper subset of $\partial B$ and that the quotient  
$(\partial B\sm \Lambda)/\Gamma$ has a compact component $A$.  
Let $E$ be the end of $M$ corresponding to $A$ and assume  
that $M\sm E$ has finite volume. Then $A$ is embeddable in some $\C^N$ and 
$M$ can be compactified to a strongly pseudoconvex domain in a projective variety 
by adding an exceptional analytic set. 
\end{cor} 
In \cite{Na-R} the case $n\geqslant3$ is considered and it is shown that the volume of $M\sm E$ is automatically finite. 
Indeed, if $n\geqslant3$, $\partial M$ is embeddable by \cite{Bo} and the method described above, combined with 
the holomorphic Morse inequalities of \cite{N-T}, show that  
$\operatorname{vol}(M\sm E)<\infty$.
Napier and Ramachandran also remark \cite[Remark, p.\,392]{Na-R} that one could prove an analogue of Corollary \ref{Burnsintro}
in the case $n\geqslant3$, by using the techniques of Siu-Yau \cite{SY:82}.


Our second application of Theorem \ref{embeddingintro} is an embedding  
result for Sasakian manifolds (the definition of this particular class  
of strongly pseudoconvex manifolds will be given in \S 5).  
Let $X$ be a Sasakian manifold. In \cite{Bi-H}, O. Biquard and M. Herzlich  
construct a K\"ahler metric on the product $(0,\infty )\times X$, called  
asymptotically complex hyperbolic. $X$ can then be viewed as the  
pseudoconvex boundary of some K\"ahler manifold, and therefore our technique  
can be applied to prove  
\begin{theoreme}\label{Sasakianintro} 
Let $X$ be a Sasakian manifold of dimension at least $3$.  
Then there is a CR embedding of $X$ in $\C ^N$ for some integer $N$. 
\end{theoreme} 
This paper is organized as follows. In \S 2, we disscuss the 
Hodge decomposition for the space of square integrable forms. In \S 3-4,  
we prove Theorem \ref{embeddingintro} and \ref{pinched}.  
Finally in \S 5, we review some basic facts about Sasakian manifolds,  
explain the construction of asymptotically hyperbolic manifolds, and prove  
Theorem \ref{Sasakianintro}. 

\section{Hodge decomposition} 
Let $(M,g)$ be a complex complete Hermitian manifold. 
For any integers $p$ and $q$, we denote by $C^{p,q}_0(M)$ the space of smooth 
$(p,q)$-differential forms with compact support, and by $L_2^{p,q}(M)$ the 
space of square integrable $(p,q)$-forms. Now consider the operator 
$\overline{\partial}+\overline{\partial}^*$, where $\overline{\partial}$ acts 
as an unbounded operator on square integrable forms, and  
$\overline{\partial}^*$ is the Hilbert space adjoint of  
$\overline{\partial}$. The $L_2$-kernel $\harm ^{p,q}(M,g)$ of  
$\overline{\partial}+\overline{\partial}^*$ is by definition the space of  
$L_2$ harmonic $(p,q)-$forms. 
In order to solve the $\overline{\partial}$-equation, one is primarily 
interested to know whether $\harm ^{p,q}(M,g)$ is finite dimensional and 
$\overline{\partial}$ has a closed range as an unbounded operator acting on 
$L_2$-forms. For example, this is the case if $\overline{\partial} 
+\overline{\partial}^*$ (or equivalently $\overline{\square}  
:=(\overline{\partial} 
+\overline{\partial}^*)^2$) is Fredholm on its domain. It is well known that  
we he have the following characterization of Fredholmness (see for example 
\cite{An}).  
\begin{prop}\label{Ancona} 
The operator $(\overline{\partial} +\overline{\partial}^*)^2$ is  
Fredholm in bidegree $(p,q)$ if and only if we have a Poincar\'e inequality  
at infinity: there  
exist a compact subset $K$ of $M$ and a constant $C>0$ such that 
\begin{equation}\label{Poincare} 
\forall u\in C^{p,q}_0(M\sm K),\, C||u||^2_{L_2}\leqslant ||\overline{\partial} 
u||^2_{L_2}+||\overline{\partial}^* u||^2_{L_2}. 
\end{equation} 
\end{prop} 
Our discussion applies to weighted spaces too. Namely, let $\xi :M\to \R$ be  
a smooth function defined on $M$, and let $L_{2,\xi}(M)$ be the Hilbert  
space completion of $C^{*,*}_0(M)$ with respect to the weighted norm 
$$\Vert u \Vert _{L_{2,\xi}}^2=\int _M \vert u\vert ^2 e^{-\xi}\,  
\rm{d}v_g,$$ 
where $\rm{d}v_g$ is the Riemannian measure.  
The operator $\overline \partial$ acts as an unbounded operator on  
$L_{2,\xi}(M)$, and its adjoint is denoted by $\overline{\partial }^*_\xi$;  
we have 
$$\overline{\partial }^*_\xi=e^\xi\overline{\partial }^*e^{-\xi}.$$ 
The space of $(\overline \partial+\overline{\partial }^*_\xi)$-harmonic $(p,q)$ forms is 
$$\harm ^{p,q}_{\xi}(M,g)=\{ u\in L_{2,\xi}^{p,q}(M):\,\overline{\partial} 
u=\overline{\partial}^*_\xi u=0\} .$$  
 
\noindent 
A key property of Fredholmness is the following:
\begin{prop}\label{NPI} 
Let $(M,g)$ be a complete Hermitian complex manifold. Assume that 
$\overline{\partial}+\overline{\partial}^*_\xi$ is Fredholm 
in bidegree $(p,q)$. Then $\harm ^{p,q}_\xi(M,g)$ is finite dimensional,  
and we have the following strong orthogonal Hodge decomposition 
\begin{equation}\label{hodge} 
L_{2,\xi}^{p,q}=\harm ^{p,q}_\xi(M,g)\oplus \im(\overline{\partial} :  
L_{2,\xi}^{p,q-1}\to L_{2,\xi}^{p,q})\oplus  
\im(\overline{\partial} ^*_\xi:  
L_{2,\xi}^{p,q+1}\to L_{2,\xi}^{p,q})\,. 
\end{equation} 
\end{prop} 
 
For the proof of Theorem \ref{pinched} we need the following example of manifolds  
for which the operator 
$(\overline \partial+\overline{\partial }^*)^2$ is Fredholm  
at certain degrees.  
\begin{prop}[Gromov]\label{Gromov} 
If $(M,g)$ is a complete K\"ahler manifold whose fundamental form is  
$d$(bounded) (at least outside a compact subset) in the sense of  
Gromov, the operator $(\overline{\partial}+\overline{\partial}^*)^2$ is Fredholm
in bidegree $(p,q)$, if $p+q\neq\dim_{\C}(M)$. 
\end{prop} 
Recall that the hypothesis means that for a bounded $1-$form $\eta$,  
the fundamental form $\omega$ of the metric $g$ satisfies $\omega =d\eta$  
outside a compact subset $K$. Then the method of proof of  
\cite[Theorem 1.4.A]{Gr}, using the Lefschetz decomposition theorem,  
shows that for some constant $C>0$ depending only on the dimension of $M$  
and the bound on the norm of $\eta$, we have 
$$\forall p+q\neq \dim_{\C}(M),\, \forall u \in C^{p,q}_0(M\sm K),\,  
C||u ||^2_{L_2}\leqslant \lan \Delta u ,u \ran _{L_2} =2\lan \overline{\square}u , 
u \ran _{L_2}.$$ 
Here, the last equality is a consequence of the K\"ahler assumption, 
for we then have $\Delta =2\overline{\square}$, and after integrating by  
parts, we get the estimates of \eqref{Poincare}.  

\section{Embedding of some strongly pseudoconvex CR manifolds} 
The aim of this section is to prove Theorem \ref{embeddingintro}. We first review
the notion of CR structure and complex manifold with strongly pseudoconvex boundary.
In Lemma \ref{metric}, we construct a complete metric and a suitable weight $\xi$ on a strongly pseudoconvex subdomain $\Omega_0$ of $M$, such that the operator $(\overline{\partial}+\overline{\partial}^*_\xi)^2$ is Fredholm in bidegree $(0,1)$. Then Lemma \ref{pic} provides peak holomorphic functions at $\partial\Omega_0$. We use these peak functions in the proof of Theorem \ref{embeddingintro} in order to fill $\partial\Omega_0$ and then $\partial M$.

Let $\Omega$ be a strongly pseudoconvex domain in a complex manifold $M$. From the  
complex structure of $M$, we can build on the boundary of $\Omega$  
a partial complex structure which is called a Cauchy-Riemann or CR structure. More generally, 
let $Y$ be a smooth 
orientable manifold of (real) dimension $(2n-1)$. A \emph{CR structure} on $Y$ is 
an $(n-1)$-dimensional complex subbundle ${H_{(1,0)}Y}$ of the complexified tangent 
bundle $T_\C Y$ such that  
$$H_{(1,0)}Y\cap \overline{H_{(1,0)}Y}=\{ 0\},$$ 
and such that $H_{(1,0)}$ is integrable as a complex subbundle of $T_\C Y$  
(i.e. if $u$ and $v$ are sections of $H_{(1,0)}Y$, the Lie bracket $[u,v]$  
is still a section of $H_{(1,0)}Y$).  
 
If $Y$ is a CR manifold, then its Levi distribution $H$ is the real subbundle 
of $TY$ defined by $H=Re\{ H_{(1,0)}Y\oplus \overline{H_{(1,0)}Y}\}$.  
There exists on $H$ a complex structure $J$ given by $J(u+\overline{u})= 
\sqrt{-1}(u-\overline{u})$, with $u\in H_{(1,0)}Y$.  
As $Y$ is orientable, the real line bundle $H^{\perp} \subset T^*Y$ admits a  
global nonvanishing section $\theta$. The CR structure is said to be  
\emph{strongly pseudoconvex} if $d\theta (.,J.)$ defines a positive definite  
metric on $H$. Notice that in this case,  
$\theta \wedge (d \theta)^{n-1}\ne 0$, and $\theta$ defines a real contact  
structure on $Y$.  

We need the notion of complex manifold with strongly pseudoconvex boundary.
\begin{definition}\label{bdy}
A complex manifold $M$ with strongly pseudoconvex boundary is a real  
manifold with boundary, of dimension $2n$,
satisfying the following conditions: (i) the interior 
$\operatorname{Int}M=M\smallsetminus\partial M$ has an integrable 
complex structure and 
(ii) for each point $x\in\partial M$ there exist a neighborhood $U$ in $M$,
a strongly pseudoconvex domain $D\subset\C^n$ with smooth boundary, and a diffeomorphism
$h$ from $U$ onto a relatively open subset $h(U)$ such that $h(\partial U)\subset \partial D$
and $h$ is biholomorphic from $\operatorname{Int}U$ to $\operatorname{Int}h(U)$. 
\end{definition}  
From this definition we infer:
\begin{consequence}\label{psh}
The complex structure induces an integrable Cauchy-Riemann 
structure on the boundary $\partial M$. Moreover, if $\partial M$ is compact, 
there exists a defining function
$\varphi:M\to(-\infty,c]$ such that $\partial M=\{\varphi=c\}$, with the properties:
(1) its Levi form is positive definite on the holomorphic tangent space of $\partial M$ and  
(2) $\varphi$ is strictly plurisubharmonic on $\{c_0<\varphi<c\}$.
We can assume that $c_0<0<c$ and that $0$ is a regular value of $\varphi$ and consider 
the strongly pseudoconvex domain $\Omega _0=\{ \phi <0 \}$. 
\end{consequence}

For the lemmas leading to Theorem \ref{embeddingintro}, we construct    
a smooth function $\lambda:(-\infty,0)\to\R_+$ such that $\lambda=0$ in a neighborhood of
$-\infty$ and $\lambda^\prime>0$, $\lambda^{\prime\prime}>0$ in a neighborhood of $0$. We will impose two conditions
on the growth of $\lambda$ near $0$. 
Let $0<\delta<C$ be two constants and let $\psi:\Omega _0 \to\R$ be a smooth function
such that 
$$\psi = \left\{ \begin{array}{ll} 
     -(n+2)\log (-\phi) & \textrm{on $\{\phi>-\delta\}$\,,}\\ 
     0 & \textrm{on $\{\phi<-C\}$\,.} 
     \end{array} \right.$$ 
First, we may choose $\lambda$ such that
\begin{equation}\label{cond1}
\lambda(\phi)\geqslant \psi\quad\text{on $\Omega_0$}
\end{equation}

\begin{lemma}\label{metric} 
Assume that $M$ can be endowed with a complete  
Hermitian metric such that the operator $(\overline{\partial}+ 
\overline{\partial}^*)^2$ is Fredholm in bidegree $(0,1)$. Then there exists on  
$\Omega _0$ a complete Hermitian metric 
(sending the boundary at infinity) for which the weight $e^{-\lambda(\phi)}$  
is integrable near the boundary and the operator 
$(\overline \partial+\overline \partial ^*_{\lambda(\phi)})^2$ is  
Fredholm in bidegree $(0,1)$. 
\end{lemma} 
\begin{proof} 
Consider the complete metric $g$ on $\Omega _0$ obtained by gluing the  
original metric of $M$ with the K\"ahler metric defined near the boundary  
$\partial \Omega _0$ by the K\"ahler form 
\begin{eqnarray*} 
\omega &=& -\sqrt{-1}\;\partial \overline{\partial}\log (-\phi ),\\ 
       &=& -\sqrt{-1}\;\frac{\partial \overline{\partial}\phi}{\phi}+ 
\sqrt{-1}\;\frac{\partial \phi \wedge \overline{\partial}\phi}{\phi  ^2}. 
\end{eqnarray*} 
Strict plurisubharmonicity of $\phi$ ensures that $\omega$ is positive. 
Near the boundary, 
the volume of $g$ behaves like $O((-\phi) ^{-(n+1)})$ (see e.g. \cite[Lemma 4.2]{Y}), so that $e^{-\psi}=(-\phi)^{n+2}$ is integrable around 
$\phi=0$. 
By \eqref{cond1}, $e^{-\lambda(\phi)}\leqslant e^{-\psi}$, hence $e^{-\lambda(\phi)}$ is also integrable around 
$\phi=0$. 

We now prove that the operator  
$(\overline \partial+\overline \partial ^*_{\lambda(\phi)})^2$  
associated to this metric is Fredholm in 
bidegree $(0,1)$. We have to show that the estimates \eqref{Poincare} hold for  
smooth $(0,1)$-forms  which are compactly supported at infinity.
 
If we are sufficiently far from the boundary $\partial \Omega _0$,  
then the new metric is equal to the original metric of $M$, and the weight  
$e^{-\lambda(\phi)}$ is identically one. 
Therefore, by assumption and by Proposition \ref{Ancona}, these estimates are valid for forms with compact support in $\Omega _0$, sufficiently far  
from the boundary $\partial \Omega _0$.  
 
It remains to prove \eqref{Poincare} for a $(0,1)$-form $u$ with compact  
support near the boundary, i.e.
in $\{\phi>-\delta\}$. 
We shall use the Bochner-Kodaira formula.
We endow the trivial line bundle $E=\Omega_0\times\C$ with the Hermitian  
metric $h=e^{-\lambda(\phi)}$.
Consider the holomorphic vector bundle $\widetilde{E}=E\otimes\Lambda^nTM$.
We denote by $\sim\,:\Lambda^{0,1}T^*M\otimes E\to\Lambda^{n,1}T^*M\otimes\widetilde{E}$, $u\mapsto
\widetilde{u}$ the natural isometry which sends a $(0,1)$-form with values in $E$ to a
$(n,1)$-form with values in $\widetilde E$.

 To avoid excessive subscripts, we will drop $\lambda (\phi )$ from the notation, 
and write e.g. $\overline \partial ^*$ instead of $\overline \partial ^*_{\lambda (\phi )}$. 
For a $(0,1)$-form $u$ with compact support, we have 
\begin{equation}\label{tilde}
\norm{\overline\partial u}^2+\norm{\overline\partial^* u}^2=\norm{\overline\partial\widetilde{u}}^2+
\norm{\overline\partial^*\widetilde{u} }^2
\end{equation}
By the Bochner-Kodaira formula for the $(n,1)$-form $\widetilde{u}$ we have
\begin{equation}\label{bk}
\norm{\overline\partial\widetilde{u}}^2+\norm{\overline\partial^*\widetilde{u} }^2\geqslant
\left(\big[\sqrt{-1}\Theta(\widetilde{E}),\Lambda\big]\widetilde{u},\widetilde{u}\right)
\end{equation}
where 
$$\sqrt{-1}\Theta(\widetilde{E})=\sqrt{-1}\Theta(E)+\operatorname{Ric}\omega=
\sqrt{-1}\partial \overline \partial\lambda(\phi)+\operatorname{Ric}\omega$$ 
is the curvature of $\widetilde{E}$, $\operatorname{Ric}\omega$ is the Ricci curvature 
and $\Lambda$ is the interior product with $\omega$.
Since $\widetilde{u}$ is of bidegree $(n,1)$, relations \eqref{tilde} and \eqref{bk} imply
\begin{equation}\label{bk1}
\norm{\overline\partial u}^2+\norm{\overline\partial^* u}^2\geqslant
\left(\sqrt{-1}\Theta(\widetilde{E})\wedge\Lambda\widetilde{u},\widetilde{u}\right)
=\left(\big(\sqrt{-1}\partial \overline \partial\lambda(\phi)+\operatorname{Ric}\omega\big)\wedge
\Lambda\widetilde{u},\widetilde{u}\right)
\end{equation}
Since $\partial \overline \partial\lambda(\phi)=\lambda^\prime(\phi)\partial \overline{\partial}\phi+
\lambda^{\prime\prime}(\phi)\partial\phi\wedge\overline{\partial}\phi$ we can choose $\lambda$ increasing 
enough so that
\begin{equation}\label{bk2}
\sqrt{-1}\partial \overline \partial\lambda(\phi)+\operatorname{Ric}\omega\geqslant\omega
\end{equation}
in a neighborhood $\{\phi>-\delta\}$ of $\partial\Omega_0$.
Estimates \eqref{bk1} and \eqref{bk2} entail
\begin{equation}\label{bk3}
\norm{\overline\partial u}^2+\norm{\overline\partial^* u}^2\geqslant
\left(\omega\wedge\Lambda\widetilde{u},\widetilde{u}\right)\geqslant\norm{\widetilde{u}}^2=\norm{u}^2\,.
\end{equation}
for $u$ supported in $\{\phi>-\delta\}$.  
This finishes the proof of the lemma. 
\end{proof}  
 
\noindent
The following result is the analog of \cite[Corollary 2.6]{M-D}. 
\begin{lemma}\label{pic}  
Assume that $M$ can be endowed with a complete 
Hermitian metric such that the operator $(\overline{\partial}+ 
\overline{\partial}^*)^2$ is Fredholm in bidegree $(0,1)$.  
Let $p$ be a point of the boundary $\partial \Omega _0$, and $f$ a  
holomorphic function defined in a neighborhood of $p$ such that  
$\{f=0 \}\cap \overline{\Omega} _0=\{ p\}$. Then for any $m$ big enough,  
there exist a function $g\in \mathcal{O}(\Omega _0 )\cap C^\infty  
(\overline{\Omega _0}\sm \{ p\})$, a smooth function $\Phi$ on a  
neighborhood $V$ of $p$, and constants $a_1,\ldots ,a_{m-1}$ such that 
$$g=\frac{1+a_{m-1}f+\ldots +a_1f^{m-1}}{f^m}+\Phi$$ 
on $V\cap \Omega _0$. In particular, we have $\lim _{z\to p}|g(z)|=\infty$. 
\end{lemma} 
\begin{proof} We follow the same line of reasoning as in  
\cite[Corollary 2.6]{M-D} with only minor changes. For $\epsilon >0$, set 
$$\Omega _\epsilon =\{\phi <\epsilon\},$$ 
and choose $\epsilon$ small enough so that $\Omega _\epsilon$ is  
strongly pseudoconvex. Replacing $\Omega _0$ by $\Omega _\epsilon$ in  
Lemma \ref{metric}, we can find a complete metric $g$ on $\Omega _\epsilon$  
and a suitable function $\xi :\Omega _\epsilon \to \R$, of the form  
$\xi=\lambda(\phi-\epsilon)$, 
with the properties stated in that lemma. 
By Proposition \ref{NPI} we have  
\begin{equation}\label{dec} 
\ke(\overline{\partial}:L^{0,1}_{2,\xi}(\Omega _\epsilon)\to L^{0,2}_{2,\xi} 
(\Omega _\epsilon)) =  
\harm ^{0,1}_\xi(\Omega _\epsilon,g) \oplus \im(\overline{\partial} :  
L^{0,0}_{2,\xi}(\Omega _\epsilon)\to L^{0,1}_{2,\xi}(\Omega _\epsilon)), 
\end{equation} 
and  
\begin{equation}\label{dim} 
\rm{dim}(\harm ^{0,1}_\xi(\Omega _\epsilon,g))<\infty . 
\end{equation} 
Now, let $U$ be a small neighborhood of $p$ where $f$ is defined, and let  
$\psi$ be a smooth function with compact support in $U$, such that $\psi =1$  
in a neighborhood $V$ of $p$. Set 
$$h_m=\psi /f^m \quad \textrm{on $U$ and $0$ on $M\sm U$}$$ 
and 
$$v_m=0 \quad \textrm{on $V$ and $\overline \partial h_m $ on $M\sm V$}.$$ 
If $\epsilon$ is sufficiently small, $v_m$ is smooth on  
$\overline{\Omega _\epsilon}$, with compact support  
(not disjoint from the boundary). In particular, 
as the weight $e^{-\xi}$ is assumed to be integrable near  
$\partial \Omega _\epsilon$, $v_m$ is in $L^{0,1}_{2,\xi}(\Omega _\epsilon)$. 
Moreover, we have $\overline \partial v_m=0$ on $\Omega _\epsilon$.  
By \eqref{dec} and \eqref{dim}, the codimension of the image  
$\im(\overline{\partial} :  L^{0,0}_{2,\xi}\to L^{0,1}_{2,\xi})$ in the kernel  
$\ke(\overline{\partial}:L^{0,1}_{2,\xi} \to L^{0,2}_{2,\xi})$ is finite,  
so that for every $m$ big enough, there are constants $a_1,\ldots ,a_{m-1}$  
such that $v=v_m+a_{m-1}v_{m-1}+\ldots +a_1v_1$ belongs to this image.  
Hence there is a function $\Phi '$ such that $\overline \partial \Phi '=-v$,  
and by elliptic regularity, $\Phi '$ is smooth on $\Omega _\epsilon$.  
Set $h=h_m+a_{m-1}h_{m-1}+\ldots +a_1h_1$ and $g=h+\Phi '$.  
Then $\overline \partial g=0$ on $\Omega _\epsilon \sm \{ f=0\}$,  
$g\in\mathcal{O}(\Omega _0 )\cap C^\infty (\overline{\Omega _0}\sm \{ p\})$,  
and the function $\Phi$ in the lemma is equal to $\Phi '$ on $V$. 
\end{proof}  
\noindent 
We can now prove Theorem \ref{embeddingintro}. 
\begin{proof}[Proof of Theorem \ref{embeddingintro}]  
Let us consider a strictly plurisububharmonic function $\varphi$
as in Consequence \ref{psh}.
\comment{
As before we can assume that $M$ is a domain in larger  
differentiable manifold $M'$ and for some smooth function  
$\phi:M'\to\R$ and some $c>0$, 
we have $M=\{x\in M':\phi(x)\leqslant c\}$, $\partial M= 
\{x\in M':\phi(x)=c\}$. Since for $d<c$ close to $c$ the domains 
$\{x\in M':\phi(x)\leqslant d\}$ are strongly pseudoconvex, we may assume  
that there exists a domain $\Omega_0\subset 
\operatorname{Int}M$ with $\Omega _0=\{ \phi <0 \}$ and $\phi$ is a  
strongly plurisubharmonic function in the neignbourhood of $\partial\Omega_0$.}
Let $\epsilon\in(0,c)$ be sufficiently small. By Lemma \ref{pic}, we can construct 
peak holomorphic functions for each point of  
$\partial \Omega _{\epsilon}=\{ \phi =\epsilon \}$.  
Using this fact, it is shown in \cite[Proposition 3.1]{M-D} that if  
$\delta >0$ is sufficiently small, the holomorphic functions on  
$\{ \phi <\epsilon \}$ separate points and give local coordinates on  
$\{\epsilon-\delta<\phi <\epsilon\}$.  
By the method of \cite[Proposition 3.2]{AS} (see also \cite[\S3]{M-D}) 
there exists a compact Stein space $S$ with boundary  
(and with at most isolated singularities)  
and an embedding of $\{\epsilon-\delta<\phi <\epsilon\}$ as an open set 
in $S$ such that $\partial S=\{\phi=\epsilon\}$.  
By gluing  
$\{\epsilon-\delta<\phi<c\}$ with $S$ along 
$\{\epsilon-\delta<\phi<\epsilon\}$ we obtain a Stein space 
$S'$ with boundary $\{\phi=c\}$. Using \cite[Theorem 0.2]{Heu} (see also Ohsawa \cite{Oh}) 
the space $S'$ can be embedded as a domain with boundary in a larger 
Stein space $S''$ such that $\partial S'$ is a hypersurface in $S''$. 
By the embedding theorem of Remmert-Bishop-Narasimhan \cite{Na:60}, $S''$ admits a proper  
holomorphic embedding in $\C^N$ for some $N$. Restricting this embedding 
to $\partial S'=\{\phi=c\}=\partial M$  
we obtain the conclusion of the theorem. 
\end{proof}  
\comment{  
\begin{remark}
From \cite[Proposition 4.4]{Na-R} it follows that $\partial M$ is connected and there is one end 
corresponding to $\partial M$. 
\end{remark}
 
\begin{remark} 
It is not necessary to assume that $\Omega _0$ is a domain in a larger  
manifold. We can instead begin directly with a complex manifold $\Omega _0$  
with a compact strictly pseudoconvex boundary. If we assume that it has a  
complete K\"ahler (outside a compact subset) metric for which the operator  
$\overline{\partial}+\overline{\partial}^*$ is Fredholm in  
bidegree $(0,1)$, we can still embed $\partial \Omega _0$ 
as the preceding proof still applies. 
\end{remark} 
}
 
\section{Compactification of manifolds with pinched negative curvature} 
 
Let us first of all review some standard facts about manifolds  
with pinched negative curvature  
(see the works Eberlein \cite{E} and  
Heintze-Im Hof \cite{H-I}, or the book \cite{B-G-S}). Thus, let $(X,g)$ be any complete manifold with  
pinched negative sectional curvature $\sigma$, that is, we assume that  
there are positive constants $a$ and $b$ such that  
$$-b^2\leqslant \sigma \leqslant -a^2<0.$$  
For a compact subset $K$ of $X$, an unbounded connected component of  
$X\sm K$ is called an end of $X$ (with respect to $K$). If $K_1\subset K_2$  
are two compact subsets, the number of ends with respect to $K_1$ is at  
most the number of ends with respect to $K_2$, so that we can define the  
number of ends of $X$. Namely, $X$ is said to have finitely many ends if  
for some integer $k$, and for any $K\subset X$, the number of ends with  
respect $K$ is at most $k$. The smallest such $k$ is called the number of  
ends of $X$, and then there exists $K_0\subset X$ such that the number  
of ends with respect to $K_0$ is precisely the number of ends of $X$.  
If no such $k$ exists, we say that $X$ has infinitely many ends. 

Let $E$ be an end of $X$, and assume that $E$ has finite volume.  
$E$ is called a cusp and is diffeomorphic to a product  
$(0,\infty )\times \Sigma$ where $\Sigma$ is a connected closed manifold.  
The slices $\{ t\} \times \Sigma$ are the level sets of a Busemann function  
$r$ associated to $E$. Moreover, the metric restricted to $E$ is given by  
$$g\vert _{E}=dr ^2+h_r,$$ where $h_r$ is a family of metrics on $\Sigma$.  
Finally, standard Jacobi field estimates show that for all $r$,  
we have $e^{-br}h_0\leqslant h_r \leqslant e^{-ar}h_0$.  
 
\begin{proof}[Proof of Theorem \ref{pinched}]  
For the proof of 1), we shall prove that the hypothesis of Theorem \ref{embeddingintro} are fullfiled.  
 
Let $U$ be a neigbourhood of $\partial M$ such that $M\sm U$ has finite volume. Then $M\sm U$ has 
finitely many ends, each of them being a cusp. 

Using the Busemann functions associated to each of these  
cusps, it is shown for example in \cite[Lemma 3.2]{Y} that the K\"ahler form of the  
metric is $d$(bounded) on each cusp. 
Let us explain the details of the construction. 
Let $E$ be a cusp of $M$, and let $r$ be a Busemann function associated to  
$E$ such that $E$ is diffeomorphic to a product $(0,\infty)\times \Sigma$,  
with $\Sigma$ a closed connected manifold. Let $\phi _t$ be the flow of the  
gradient $\nabla r$ of $r$. If $x=(s,y)$ is a point of  
$(0,\infty )\times \Sigma$, we simply have 
$$\phi _t(s,y )=(s+t,y ).$$ 
The properties of the differential of this flow are well understood.  
First, for $x\in E$, we have $d_x\phi _t(\nabla r(x))=\nabla r(\phi _t(x))$.  
Let $u$ be a tangent unitary vector at $x$, and orthogonal to $\nabla r$.  
Then $t\to d_x\phi _t(u)$ is a stable Jacobi field along the geodesic  
$t\to(s+t,y)$, orthogonal to $\nabla r$, and is equal to $u$ at $t=0$.  
As the sectional curvature is less than or equal to $-a^2$, the classical  
estimates on Jacobi fields \cite{H-I} show that 
$$\vert d_x\phi _t(u)\vert \leqslant e^{-a\,t}.$$ 
Now, if $\omega$ is the K\"ahler form of the metric, it is closed, so that  
by Cartan formula, we have 
$$\phi _t^*\omega -\omega =d\int _0^t \phi _s^*(i_{\nabla r}\omega )\, ds.$$ 
By using the estimate of the differential of $\phi _t$, we can take the  
limit as $t$ goes to infinity and get that $\omega =d\eta$ on the cusp $E$,  
with 
$$\eta =-\int _0^\infty \phi _s^*(i_{\nabla r}\omega )\, ds.$$ 
The same estimate imply that $\eta$ is bounded.  

From Proposition \ref{Gromov} it follows  
that $(\overline \partial+\overline \partial ^*)^2$ is  
Fredholm in bidegree $(0,1)$.  
\noindent 
By Theorem \ref{embeddingintro} this achieves the proof of 1). 

We consider conclusion 2).  
The first observation is that $M$ can be compactified, which 
allows the use of the extension theorem of Kohn-Rossi. 
Indeed, let $E_1,\dotsc,E_m$ be the cusps of  
$M\smallsetminus U$. We fix some end $E_j$ and consider the associated Busemann 
function $r:E_j\to(0,\infty)$. It follows from 
\cite[Proposition 1]{SY:82} that $-r:E_j\to(-\infty,0)$ 
is a strictly plurisubharmonic proper function (note that for the Busemann function, Siu and Yau use the opposite sign convention).  
From Lemma \ref{pic} for $\Omega_0=\{-r<c\}$, $c<0$, 
we infer as in the proof of Theorem \ref{embeddingintro} that  
each strip $\{c<-r<0\}$ can be compactified to a Stein space 
with isolated singularities $S_{j,c}$. By the uniqueness of the Stein  
completion \cite[Corollary 3.2]{AS} all spaces $S_{j,c}$, $c<0$, 
are biholomorphic to a fixed Stein space $S_j$. 
Letting $c\to-\infty$, we see that $E_j$ is biholomorphic to 
an open set of $S_j$. Therefore all ends of $M\smallsetminus U$, and together 
with them also $M$, can be compactified. 
 
In particular, the Kohn-Rossi theorem \cite{KR:65} shows that 
every holomorphic function defined in a  
neighborhood of $\partial M$ extends to a holomorphic function on 
$M$. As a by-product we obtain that $\partial M$ is connected. 
 
The next step is to glue $M$ to a pseudoconcave projective manifold. 
Since $\partial M$ is embeddable in $\C^N$  
we can apply \cite[Theorem 1.2]{Le2} which allows to  
assume that the Stein space $S''$ constructed in the proof of Theorem  
\ref{embeddingintro} is an open set in an affine algebraic variety, hence also 
in a projective variety $X$. 
We use now the notations from the proof of Theorem \ref{embeddingintro}.
We set $N=\{\varepsilon-\delta<\varphi\leqslant c\}$  
and glue the 
manifolds $M$ and $(X\smallsetminus S')\cup N$ along $N$. The resulting manifold 
will be denoted by $\widehat{M}$. 
Hence $M$ is a domain with compact strongly 
pseudoconvex boundary in $\widehat{M}$.
 
Since $S''$ is an affine space in some $\C^N$, we can regard the embedding 
of $N'=\{\varepsilon-~\delta<\varphi<\varepsilon\}$
in $X$ as a map with values in $\C^N$. Now $M$ can be compactified to 
a compact strongly pseudoconvex domain, so the extension theorem of Kohn-Rossi, 
applied to the components of this embedding, show that the embedding
extends to a holomorphic map from $M$ to $\C^N\subset\mathbb{P}^N$. 
 
Pulling back the hyperplane line bundle of $\mathbb{P}^N$ through this map, we obtain
a line bundle $E\to\widehat{M}$ which is semi-positive on $\widehat{M}$
and positive on $(X\smallsetminus S')\cup N$. 

A partition of unity argument delivers a Hermitian metric on  
$\widehat{M}$ which agrees with the 
original metric $\omega$ of $\operatorname{Int}M$ on say $\{\varphi<\varepsilon\}$. 
With respect to this metric 
the canonical bundle of $M$ is positive on $\{\varphi<\varepsilon\}$. Hence,
the bundle $L=E^k\otimes K_{\widehat{M}}$ is positive on $\widehat{M}$ for $k$ sufficiently large.  
Moreover, the curvature $\sqrt{-1}\Theta(L)$ of $L$ dominates $\omega$ on $\{\varphi<\varepsilon\}$ 
and therefore $\Theta(L)$
is a complete metric on $\widehat{M}$. The $L_2$ estimates of  
H\"ormander-Andreotti-Vesentini \cite{De:01}
produce sections of $\oplus_{\nu}H^0(L^{\nu}\otimes K_{\widehat{M}})$  
that separate points and give local coordinates 
on $\widehat{M}$.

On the other hand the manifold $\widehat{M}$ is $1$-concave in the sense  
of Andreotti-Grauert (we use again the plurisuperharmonicity of the Busemann 
function on each cusp).
Therefore, the embedding theorem of Andreotti-Tomassini \cite{AnTo:70}  
implies that the ring
$\oplus_{\nu}H^0(L^{\nu}\otimes K_{\widehat{M}})$ gives an embedding of $\widehat{M}$
in some projective space $\mathbb{P}^N$. By \cite{An:63} the projective closure of 
$\widehat{M}$ in $\mathbb{P}^N$ is a projective variety $\widetilde{M}$ of the same dimension.

We have thus found a compactification $\widetilde{M}\subset\mathbb{P}^N$  
of $\widehat{M}$ and therefore of $M$. 
The set $\widetilde{M}\smallsetminus M$ is a pluripolar set, namely the set where
the plurisubharmonic function $-r$ takes the value $-\infty$. 
By Wu's theorem \cite{G-W}, any simply connected complete K\"ahler  
manifold of nonpositive sectional curvature is Stein.  
Hence the universal covering of $M$ is Stein.
It is then shown in \cite[\S\,4]{SY:82} (using the Schwarz-Pick Lemma of Yau) and in
\cite[Theorem 1.3]{M-D} (using Wermer's theorem) that $\widetilde{M}\smallsetminus M$
is an exceptional analytic set in the sense of Grauert  
\cite[Definition 3,\,p.\,341]{Gra:62}.

\noindent  
More precisely, let us consider the compact strongly pseudoconvex domain
$V=\widetilde{M}\smallsetminus(X\sm S')$ which compactifies $M$. By \cite[Satz 3,\,p.338]{Gra:62}
there exists a maximal nowhere discrete analytic set $A$ of $V$ and by
\cite[Satz 5,\,p.340]{Gra:62} there exists a Remmert reduction $\pi:V\to V'$,
which blows down $A$ to a discrete set of points.
Here $V'$ is a normal Stein space with at most isolated singularities.
Then $\pi(\widetilde{M}\smallsetminus M)=\pi(V\smallsetminus M)$ is a subset 
of $\operatorname{Sing}V'$ and each point of $\pi(\widetilde{M}\smallsetminus M)$
compactifies one of the ends $E_1$, \dots, $E_m$ of $M\sm U$. 
Thus $\widetilde{M}\smallsetminus M$ consists of those components of 
$A$ obtained by blowing-up points in $\pi(\widetilde{M}\smallsetminus M)$. 
\end{proof}  
\begin{proof}[Proof of Corollary \ref{Burnsintro}]   
As is well known, the limit  
set $\Lambda$ is the set of accumulation points of any orbit  
$\Gamma\cdot x$, $x\in B$, and is a closed $\Gamma$-invariant subset of the 
sphere at infinity $\partial B$. The complement $\partial B\sm \Lambda$ is  
precisely the set of points at which $\Gamma$ acts properly discontinuously,  
and the space $M\cup (\partial B\sm \Lambda)/\Gamma$ is a manifold with 
boundary $(\partial B\sm \Lambda)/\Gamma$ (see for example  
\cite[\S 10]{E-O}). $A$ is a compact subset of this boundary, hence there  
is a neighborhood $E$ of $A$ in $M$ which is diffeomorphic to the product  
$A \times (0,1)$. It follows that $E$ is an end of $M$, because $A$ is  
compact and connected. Actually, $E$ is a strongly pseudoconvex end, in the  
sense that its boundary $A$ at infinity is strictly pseudoconvex.  
Since $M=B/\Gamma$ is a complete manifold with sectional  
curvature pinched between $-4$ and $-1$, Corollary \ref{Burnsintro}
is an immediate consequence of Theorem \ref{pinched}.  
\end{proof}  

\begin{remark}
By \cite[Theorem 4.2]{Na-R} the following holds:
\begin{burns}[Burns, Napier-Ramachandran]\label{Burns} 
Let $\Gamma$ be a torsion-free discrete group of automorphisms of the  
unit ball $B$ in $\C ^n$ with $n\geqslant 3$ and let $M=B/\Gamma$.  
Assume that the limit set $\Lambda$ is a proper subset of $\partial B$  
and that the quotient $(\partial B\sm \Lambda)/\Gamma$ has a compact  
component $A$. Then $M$ has only finitely many ends, all of which,  
except for the unique end corresponding to $A$, are cusps.  
In fact, $M$ is diffeomorphic to a compact manifold with boundary and can  
be compactified.  
\end{burns}  
The proof is based on \cite[Theorem 4.1]{Na-R}
which shows that the finite volume hypothesis of Corollary \ref{Burnsintro}
is automatically satisfied in the case $n\geqslant 3$. The presence of the strongly
pseudoconvex boundary forces the volume to be finite, since $\partial M$ is then embeddable
by \cite{Bo}, having real dimension at least $5$.

If $n=2$ we have to assume the volume to be finite in order to obtain the 
embedding of the boundary.
It is interesting to ask whether Burns' theorem holds also in dimension 
$2$ or, equivalently, whether the compact strongly pseudoconvex component
of a set $(\partial B\sm \Lambda)/\Gamma$ is embeddable for all
torsion-free discrete groups of automorphisms of the  
unit ball $B$ in $\C ^2$.

\comment{
We would like to have similar conclusions for $n=2$, possibly with stronger  
assumptions on $M$. 
 
We now recall why the hypothesis $n\geqslant 3$ is important in the proof of 
Burns' Theorem (see \cite{Na-R} for more details). 

If we are  
able to show that $M\sm E$ has finite volume, then the  
thin-thick decomposition of $M$ as in \cite{B-G-S} 
(see the remark after the proof of Lemma 10.3 and Theorem 10.5 in this book) 
will finish the proof of  
Burns' theorem. Now the difficult part is to prove that $M\sm E$ has  
finite volume. The first step is to embed the pseudoconvex boundary $A$,  
and this is possible precisely because it has real dimension at least $5$.  
Then the main tool used by Napier and Ramachandran is the asymptotic $L^2$  
Riemann-Roch inequality of Nadel and Tsuji \cite{N-T}. 
 
As a conclusion, if we could prove that $A$ is embeddable, the proof of Napier 
and Ramachandran would also work for $n=2$. Unfortunately, we were not able 
to get such an embedding. However, if we \emph{assume} that the conclusion of 
Burns' Theorem holds for $n=2$, then we can embed $A$: 
\begin{theoreme}\label{plongement} 
Let $M$ be as in Burns' Theorem, with $n=2$, and let $E$ be the end of $M$ 
corresponding to $A$. Assume that $M\sm E$ has finite volume. Then $A$ is embeddable. 
\end{theoreme}
} 
\end{remark}
 
 
\section{Embedding of Sasakian $3-$manifolds} 
\subsection{Sasakian manifolds} 
In this section, we explain some well known facts about Sasakian manifolds. Let $X$ be a strictly pseudoconvex $CR$ manifold, with compactible complex structure $J$, and compactible contact form $\theta$. This allows us to define a Riemannian metric $g_\theta$ on $X$ given by 
$$g_\theta (.,.)=d\theta (.,J.)+\theta (.)\theta (.).$$ 
Let $R$ be the Reeb vector field associated to $\theta$, defined by 
$$i_R\theta =1,\,\,\, i_Rd\theta =0.$$  
Associated to the data $(X,\theta,R,J,g_\theta)$, there is a canonical connection $\nabla$ on $TX$, called the \emph{Webster connection} (see Tanaka \cite{T} and Webster \cite{W}), which is the unique affine connection on $TX$ such that 
\be 
 \item $\nabla g_\theta =0$, $\nabla J=0$, $\nabla \theta =0$. 
 \item For any $u$, $v$ in the Levi distribution $H$, the torsion $T$ of $\nabla$ satisfies  $T(u,v)=d\theta (u,v) R$ and $T(R,Ju)=JT(R,u)$ 
\ee 
In particular, the torsion of the Webster connection cannot vanish identically. However, we have the following definition 
\begin{definition}\label{Sas} 
A strictly pseudoconvex manifold is called a Sasakian manifold if the torsion of its Webster connection in the direction of the Reeb vector field vanishes, i.e. $T(R,.)=0$ with the notations above. 
\end{definition} 
Examples of compact Sasakian manifolds are the unit sphere in $\C ^n$ or the Heisenberg nilmanifold (see Urakawa \cite{U}). Sasakian $3$-manifolds were classified up to diffeomorphism by H. Geiges \cite{G}. They were further studied by F. Belgun in \cite{Be}: every Sasakian $3-$manifold is obtained as a deformation of some standard model (see \cite{Be} for more details).

\subsection{Asymptotically complex hyperbolic manifolds} 
In \cite{Bi-H}, O. Biquard and M. Herzlich consider a class of manifolds  
which are modelled on the complex unit ball, and are thus called  
asymptotically complex hyperbolic. This construction will allow us to get  
an embedding theorem for Sasakian manifolds, but let us first recall what  
an asymptotically complex hyperbolic manifold is. 
 
Let $X$ be a $(2m-1)-$dimensional compact manifold, $m\geqslant 2$. We assume that  
$X$ has a strongly pseudoconvex CR-structure. Let $\theta$ be a compatible contact form, and $J$ a compatible almost complex structure. $\gamma (.,.) :=d\theta (.,J.)$ is then a metric on the contact distribution. Following Biquard and Herzlich \cite{Bi-H}, we endow $\Omega :=(0,\infty )\times X$ with the metric 
\begin{equation}\label{biq-her} 
g=dr^2 +e^{-2r}\theta ^2 +e^{-r}\gamma . 
\end{equation} 
Actually, in \cite{Bi-H}, the authors consider the metric $dr^2 +e^{2r}\theta ^2 +e^{r}\gamma$ on $\Omega$, but the reason for our choice will become clear later. 
 
We can extend the almost complex structure to act on the whole tangent bundle as follows. Consider the Reeb vector field $R$ and define $$J\partial _r=e^{r}R,$$ 
 where $\partial _r$ is the unit vector field in the $r$ direction. $g$ is then a Hermitian metric with respect to $J$, i.e. $J$ is an $g-$isometry. The fundamental $2-$form associated to $g$ is 
\begin{equation}\label{form} 
\omega =d(e^{-r} \theta). 
\end{equation} 
Although $\omega$ is a closed form, in this general setting, $((0,\infty )\times X,g)$ is not necessarily a K\"ahler manifold, because $J$ is not necessarily an  integrable almost complex structure. Indeed, the proof of \cite[Proposition 3.1]{Bi-H} shows that $J$ is integrable if and only if the torsion of the Webster connection of $(X,\theta )$ in the direction of the Reeb vector field vanishes identically, i.e. if and only if $X$ is a Sasakian manifold (see Definition \ref{Sas}). 
 
From now on, $X$ will be a Sasakian manifold, so that $(\Omega ,g)$ is a K\"ahler manifold.  
Note that $\eta :=e^{-r}\theta$ is a bounded $1-$form, hence the K\"ahler  
form \eqref{form} associated to $g$ is $d$(bounded). In particular,  
by Proposition \ref{Gromov}, the operator  
$({\overline\partial}+ {\overline\partial}^*)^2$ on $(\Omega , g)$  
is Fredholm in bidegree $(0,1)$ if $\dim_{\R}X\geqslant3$. Moreover, we can also write 
$$\omega =-2\sqrt{-1}\partial {\overline \partial}r,$$ 
so that $\Omega $ has a stromgly pseudoconvex boundary $\{ 0\}\times X$ 
(this explains the choice of sign of the $r$ variable for $g$ in 
\eqref{biq-her}: with  
our choice, $\Omega$ has a strongly pseudoconvex boundary, whereas with the  
choice of \cite{Bi-H} $\Omega$ has a strongly pseudoconcave boundary).  
From Theorem \ref{embeddingintro} we can conclude: 
\begin{theoreme}\label{Sasakian} 
Let $X$ be a Sasakian manifold of dimension at least $3$. Then there is a CR embedding of $X$ in $\C ^N$ for some integer $N$. 
\end{theoreme} 
\Rk The function $-r :\Omega  \to (-\infty ,0)$ is a proper smooth strictly plurisubharmonic function. Therefore $\Omega$ is a hyperconcave end, and our theorem is also a direct consequence of \cite{M-D}. 
  

\providecommand{\bysame}{\leavevmode\hbox to3em{\hrulefill}\thinspace}
\providecommand{\MR}{\relax\ifhmode\unskip\space\fi MR }
\providecommand{\MRhref}[2]{%
  \href{http://www.ams.org/mathscinet-getitem?mr=#1}{#2}
}
\providecommand{\href}[2]{#2}

\end{document}